\documentclass[a4paper]{article}
\usepackage{mathrsfs}
\usepackage{cite}
\usepackage[top=1in, bottom=1in, left=1.20in, right=1.20in]{geometry}
\usepackage{latexsym, amsmath, amssymb}
\usepackage{graphicx,epsfig}

\usepackage{amsfonts}
\usepackage{amscd}
\usepackage{amsbsy}
\usepackage{bm}
\usepackage{color}

\newtheorem{theorem}{Theorem}[section]

\newtheorem{lemma}{Lemma}[section]

\numberwithin{equation}{section}

\newcommand{\Gs}{\sigma}

\newcommand{\qed}{\ \ensuremath{\square}}

\newcommand{\pf}{\medskip \noindent {\sl Proof}. ~ }

\newcommand{\eqnref}[1]{(\ref {#1})}

\newcommand{\RR}{\mathbb{R}}

\newcommand{\beq}{\begin{equation}}
\newcommand{\eeq}{\end{equation}}
\begin{document}
\title{{\itshape Explicit calculation of strong solution on linear parabolic equation}
\thanks{Fang is supported by NSF grants  No. NSFC70921001 and No. NSFC71210003. Deng and Li are supported by NSF grants No. NSFC11301040,}}
\author{Xiaoping Fang \thanks{Postdoctoral, Management Science and Engineering Postdoctoral Mobile Station, School of Business; School of Mathematics and Statistics, Central South University, Changsha, Hunan 410083, P. R. China. Email: fxpmath@csu.edu.cn} \quad
    Youjun Deng \thanks{Corresponding author. School of Mathematics and Statistics,
   Central South University, Changsha, Hunan 410083, P. R. China. Email: youjundeng@csu.edu.cn, dengyijun\_001@163.com}  \quad
    Jing Li \thanks{Department of Mathematics, Changsha University of Science and Technology, Changsha, Hunan 410004, P.R.China. Email: lijingnew@126.com}
}
\date{}
\maketitle

\begin{abstract}
In this paper, we give the existence and uniqueness of the strong solution of one dimensional
linear parabolic equation with mixed boundary conditions. The boundary conditions can be any kind
of mixed Dirichlet, Neumann and Robin boundary conditions. We use the extension method to get the
unique solution. Furthermore, the method can also be easily implemented as a numerical method.
Some simple examples are presented. 
\end{abstract}

\section{Introduction}
We consider the existence and uniqueness of the strong solution of the following problem
\begin{equation}
\label{eq:101}
\left\{
\begin{array}{ll}
 u_t = ku_{xx} + F(x,t) & (x,t) \in \Omega_T, \\
 u(x,0)=\mu_0(x) & x\in (0,l), \\
 u_x(0,t)=0 \ \  (\mbox{or} \ \ u(0,t)=0) & t \in (0,T], \\
 -k u_x(l,t)=\nu (u(l,t)-T_0(t))  & t \in (0,T],
\end{array}
\right.
\end{equation}
where $\Omega_T:=\{0<x<l, 0<t\leq T\}$, $k>0$, $\nu>0$. For the sake of convenience we assume that $F(x,t)$, $\mu_0(x)$ and $T_0(t)$ are all
analytic functions in space and time variables, and we shall see that the solution to (\ref{eq:101}) is also analytic.

The mathematical model (\ref{eq:101}) arises in various physical and engineering settings, in particular in
hydrology \cite{Bea:1972}, material sciences \cite{RHN:1987}, heat transfer \cite{Win:1997} and transport problems \cite{ZBe:1995}.
In the case of heat transfer, the term $F(x,t)$ stands for the heat source.
In the mean time the Robin boundary condition in (\ref{eq:101}) is quite physically meaningful. Suppose
that the end $x=l$ of the rod is exposed to air or other fluid with temperature $T_0(t)$, then $u(l,t)-T_0(t)$ is the
temperature difference between the rod and its surroundings at $x=l$. By the Newton's law of cooling, which states that
the rate of change of the temperature of an object is proportional to the difference between its own temperature
and the ambient temperature (i.e., the temperature of its surroundings), we get that $-ku_x(l,t)=\nu [u(l,t)-T_0(t)]$.
The constant $\nu>0$ is called the convection coefficient or heat transfer coefficient. To sum up,
the physical model (\ref{eq:101}), considered as heat transfer, is conduction of heat in a rod whose left end is insulated (or
held at constant temperature) and whose right end is exposed to convective heat transfer. In the following paper we shall
call problem (\ref{eq:101}) as {\em Neumann-Robin} boundary problem if boundary condition $u_x(0,t)=0$ is used and {\em Dirichlet-Robin} boundary problem
if $u(0,t)=0$ is used.

Not only doest the direct heat conduction problem attract lots of researchers' attentions, but also the inverse problem.
The inverse problems of the recovery of unknown sources in (\ref{eq:101}) from
overspecified boundary data or final overdetermination have attracted many mathematicians for decades since these
problems are of paramount importance in heat conduction processes (see, e.g., \cite{BJV:1985,Can:1984,IBu:1995,Isa:1991}).
The inverse problem for simultaneous determination of source terms $F(x,t)$ and $T_0(t)$ with final overspecified data $u(x,T)$
has been studied in \cite{Has:2007}.

The uniqueness of the solution to (\ref{eq:101}) is clear and easily proved. In the mean time, the existence of the weak solution of similar form of (\ref{eq:101}) but with different boundary conditions,
has been proved long time ago by Riesz representation theorem (see e.g. \cite{Kaw:1990,LOx:1998,Pao:1975,Wlo:1987,Zeg:2002}). The existence of
the strong solution has also been proved by using the maximum principle and comparison principle for the weak solution in Dirichlet
boundary problem (cf. \cite{WYW:2006}).
However, the existence of the strong solution of (\ref{eq:101}) is still not quite clear. Here we mainly consider the strong solution of (\ref{eq:101}).
More importantly, we shall give a direct method to solve the strong solution instead of using the numerical methods such as finite difference method,
etc.
We know that the solution of (\ref{eq:101}) without boundary condition would be
\begin{equation}
\label{eq:001}
u(x,t)  =  \int_{-\infty}^{\infty} \frac{1}{\sqrt{4k\pi t}}e^{-\frac{(x-y)^2}{4kt}} \mu_0(y) dy
+ \int_{0}^{t}\int_{-\infty}^{\infty} \frac{1}{\sqrt{4k\pi (t-s)}}e^{-\frac{(x-y)^2}{4k(t-s)}} F(y,s) dyds.
\end{equation}
Concerning on the initial condition there hold
\begin{eqnarray*}
u(x,0) & = & \lim_{t\rightarrow 0}\int_{-\infty}^{\infty} \frac{1}{\sqrt{4k\pi t}}e^{-\frac{(x-y)^2}{4kt}} \mu_0(y)dy=\mu_0(x) \\
u_x(x,0) & = & - \lim_{t\rightarrow 0}\int_{-\infty}^{\infty} \frac{x-y}{2kt\sqrt{4k\pi t}}e^{-\frac{(x-y)^2}{4kt}} \mu_0(y) dy= \mu_0'(x).
\end{eqnarray*}
In fact, if we introduce the well-known Dirac function for $\delta(x)=\frac{1}{\sqrt{4k\pi t}} e^{-\frac{x^2}{4kt}}$
while $t=0$, we can easily deduce the above results.

In this paper we explore the extension method to get the strong solution to (\ref{eq:101}).
The method can also be used in other boundary conditions such as both Neumann boundaries or Dirichlet boundaries, Neumann and
Dirichlet mixed boundaries, etc.

\section{Main results}
We present the main theorem in this paper.
\begin{theorem}
\label{th:201}
Suppose $\mu_0(x)$, $F(x,t)$ and $T_0(t)$ are all analytic functions.
Then the solution to (\ref{eq:101}) has a unique solution. In particular, if $F(x,t)=0$ and $T_0(t)=\sum_0^N d_i t^i$ then
\begin{enumerate}
  \item[(1)] Neumann-Robin boundary condition ($u_x(0,t)=0$). The solution $u$ to (\ref{eq:101}) has the form
\beq\label{eq:solnr}
u=\sum_{i=0}^N a_i\sum_{j=0}^i c_jC_{2i}^{2j} x^{2i-2j}(4kt)^j+\sum_{n=1}^{\infty} b_n e^{-\sigma_n^2kt}\cos(\sigma_n x) + C_T
\eeq
where $C_T=d_0-\frac{k}{\nu}\mu'(l)-\mu(l)$ and the function $\mu(x)$ is defined by $\mu(x)= \sum_{i=0}^N a_i x^{2i}$. The
coefficients $a_i, i=0,1,\ldots, N$ are determined by
\begin{equation}
\label{matrix}
\begin{bmatrix}
1 & l(l-\frac{2k}{\nu}) & \cdots& l^{2N-1}(l-\frac{2kN}{\nu}) \\
0 & 4kc_1 & \cdots & l^{2N-3}(lC_{2N}^{2}c_1-\frac{2k}{\nu}C_{2N}^3c_2)4k \\
\vdots & \vdots & \ddots & \vdots \\
0 & \cdots & \cdots & l(lC_{2N}^{2N-2}c_{N-1}-\frac{2k}{\nu}C_{2N}^{2N-1}c_N)(4k)^{N-1} \\
0 & 0 & \cdots & c_N (4k)^N
\end{bmatrix}
\begin{bmatrix}
a_0 \\
a_1 \\
\vdots \\
a_{N-1} \\
a_N \\
\end{bmatrix}
=
\begin{bmatrix}
d \\
d_1 \\
\vdots \\
d_{N-1} \\
d_N
\end{bmatrix}
\end{equation}
where $d=\frac{k}{\nu}\mu'(l)+\mu(l)$  and
$$c_j=\frac{(2j-1)(2j-3)\cdots 1}{2^j}\quad j=0,1,\ldots,N.$$
The coefficients $b_n, n=1,2,\ldots, \infty$ are
$$b_n= \frac{2\nu }{\nu l+k \sin^2(\sigma_nl)}\int_0^l (\mu_0(x)-\mu(x)-C_T) \cos(\sigma_n x) dx$$
with $\sigma_n>0, k \sigma_n \tan(\sigma_n l) = \nu$.
  \item[(2)] Dirichlet-Robin boundary condition ($u(0,t)=0$). The solution $u$ to (\ref{eq:101}) has the form
\beq\label{eq:soldr}
u(x,t) = \sum_{i=0}^N \tilde{a}_i\sum_{j=0}^i c_jC_{2i+1}^{2j} x^{2i-2j+1}(4kt)^j+ \sum_{n=1}^{\infty} \tilde{b}_n e^{-\sigma_n^2kt}\sin(\tilde{\sigma}_n x)+\tilde{C}_T
\eeq
where $\tilde{C}_T=d_0-\frac{k}{\nu}\tilde{\mu}'(l)-\tilde{\mu}(l)$. The function $\tilde{\mu}(x)$ is defined by $\tilde{\mu}(x)= \sum_{i=0}^N \tilde{a}_i x^{2i+1}$ and $\tilde{a}_i, i=1,2, \ldots,N$ are
determined by
\begin{equation}
\label{matrix3}
\begin{bmatrix}
l-\frac{k}{\nu} & l^2(l-\frac{3k}{\nu}) & \cdots& l^{2N}(l-\frac{k(2N+1)}{\nu}) \\
0 & (lC_3^2c_1-\frac{2k}{\nu}c_2)4k & \cdots & l^{2N-2}(lC_{2N+1}^{2}c_1-\frac{2k}{\nu}C_{2N+1}^3c_2)4k \\
\vdots & \vdots & \ddots & \vdots \\
0 & \cdots & \cdots & l(lC_{2N+1}^{2N-2}c_{N-1}-\frac{2k}{\nu}C_{2N+1}^{2N-1}c_N)(4k)^{N-1} \\
0 & 0 & \cdots & (lC_{2N+1}^{2N}c_{N-1}-\frac{2k}{\nu})(4k)^{N}
\end{bmatrix}
\begin{bmatrix}
\tilde{a}_0 \\
\tilde{a}_1 \\
\vdots \\
\tilde{a}_{N-1} \\
\tilde{a}_N \\
\end{bmatrix}
=
\begin{bmatrix}
\tilde{d} \\
d_1 \\
\vdots \\
d_{N-1} \\
d_N
\end{bmatrix}
\end{equation}
where $\tilde{d}=\frac{k}{\nu}\tilde{\mu}'(l)+\tilde{\mu}(l)$. The coefficients $\tilde{b}_n, n=1,2,\ldots,\infty$ are given by
$$\tilde{b}_n=\frac{2\nu }{\nu l+k \cos^2(\tilde{\sigma}_nl)}\int_0^l (\mu_0(x)-\tilde{\mu}(x)-\tilde{C}_T) \sin(\tilde{\sigma}_n x) dx$$
where $\sigma_n>0$ satisfy $k \tilde{\sigma}_n  = -\nu\tan(\tilde{\sigma}_n l)$.
\end{enumerate}
\end{theorem}
In what follows, we shall prove the Theorem by extension method introduced in this paper. We shall consider the
Neumann-Robin problem and Dirichlet-Robin problem, separately.
\subsection{Neumann-Robin boundary condition}
For the existence of the strong solution of (\ref{eq:101}) with Neumann-Robin boundary condition we first consider
the following problem
\begin{equation}
\left\{
\begin{array}{ll}
 u_t = ku_{xx} & (x,t) \in \Omega_T, \\
 u(x,0)=\mu(x) & x\in (0,l), \\
 u_x(0,t)=0, \ \ -ku_x(l,t)=\nu(u(l,t)-T_1(t))  & t \in (0,T],
\end{array}
\right.
\label{eq:103}
\end{equation}
where $T_1(t)=T_0(t)-T_0(0)+\frac{k}{\nu}\mu'(l)+\mu(l)$ and $\mu(x)$ is an analytic function to be constructed.
We mention that in the construction of $T_1(t)$, compatibility condition on $(x=l,t=0)$ is satisfied automatically.
We consider the solution in the following form
\begin{equation}\label{eq:a01}
u(x,t)  =  \int_{-\infty}^{\infty} \frac{1}{\sqrt{4k\pi t}}e^{-\frac{(x-y)^2}{4kt}} \mu(y) dy
\end{equation}
then we have
\begin{equation} \label{eq:a02}
u_x(x,t)  =  -\int_{-\infty}^{\infty} \frac{x-y}{2kt\sqrt{4k\pi t}}e^{-\frac{(x-y)^2}{4kt}} \mu(y) dy.
\end{equation}
The strategy we use in this paper is to construct the function $\mu(x)$, $x\in(0,l)$ and then make an extension to $\RR$
to cope with the boundary conditions. For the boundary condition $u_x(0,t)=0$, we extend $\mu(x)$ to be an even function in $\RR$,
that is $\mu(x)=\mu(-x)$. Clearly we have
$$u_x(0,t) = \frac{1}{2kt\sqrt{4k \pi t}}\int_{-\infty}^{\infty} xe^{-\frac{x^2}{4kt}} \mu(x) dx =0.$$
Concerning with the boundary condition at $x=l$ we define
\begin{equation}
\label{eq:104}
g_1(t) = \frac{1}{\sqrt{4k \pi t}} \int_{-\infty}^{\infty} \left (1 - \frac{l-y}{2\nu t} \right)e^{-\frac{(l-y)^2}{4kt}} \mu(y) dy
\end{equation}
We note that $g_1(t)\in C^{\infty}((0,T))$. We suppose that $\mu(x)$ has the following form in $\RR$
\begin{equation}
\label{eq:105}
\mu(x)= \sum_{i=0}^N a_i x^{2i}.
\end{equation}
Substituting (\ref{eq:105}) into (\ref{eq:104}) we have
\begin{eqnarray*}
g_1(t) & = & \frac{1}{\sqrt{4k \pi t}} \int_{-\infty}^{\infty} \left (1 - \frac{l-y}{2\nu t} \right)
e^{-\frac{(l-y)^2}{4kt}} \sum_{i=0}^N a_i y^{2i} dy \\
& = & \frac{1}{\sqrt{4k \pi t}} \sum_{i=0}^N a_i\int_{-\infty}^{\infty} \left (1 - \frac{y}{2\nu t} \right)
e^{-\frac{y^2}{4kt}}  (l-y)^{2i}  dy \\
& = & \frac{1}{\sqrt{\pi}} \sum_{i=0}^N a_i\int_{-\infty}^{\infty} e^{-y^2}  (l-\sqrt{4kt}y)^{2i} dy
-  \frac{\sqrt{4kt}}{2\nu t\sqrt{\pi}} \sum_{i=0}^N a_i\int_{-\infty}^{\infty} y e^{-y^2}  (l-\sqrt{4kt}y)^{2i} dy ,
\end{eqnarray*}
By elementary calculation we obtain
\begin{eqnarray*}
  \int_{-\infty}^{\infty} e^{-y^2}  (l-\sqrt{4kt}y)^{2i} dy
& = & \sum_{j=0}^{i}C_{2i}^{2j}l^{2i-2j}\frac{(2j-1)(2j-3)\cdots 1}{2^j} \sqrt{\pi}(4kt)^j
\end{eqnarray*}
and
\begin{eqnarray*}
 \int_{-\infty}^{\infty} y e^{-y^2}  (l-\sqrt{4kt}y)^{2i} dy
& = & -\sum_{j=1}^{i}C_{2i}^{2j-1}l^{2i-2j+1}\frac{(2j-1)(2j-3)\cdots 1}{2^j} \sqrt{\pi}(4kt)^{j-\frac{1}{2}}.
\end{eqnarray*}
Thus by setting $c_j=\frac{(2j-1)(2j-3)\cdots 1}{2^j} $, with $c_0=1$ we have
\begin{eqnarray*}
g_1(t)& = & \sum_{i=0}^N a_i l^{2i}+ \sum_{i=1}^N a_i\sum_{j=1}^{i}C_{2i}^{2j}l^{2i-2j}c_j (4kt)^j
- \frac{1}{2\nu t} \sum_{i=1}^N a_i\sum_{j=1}^{i}C_{2i}^{2j-1}l^{2i-2j+1}c_j (4kt)^{j} \\
& = & a_0 + \sum_{i=1}^N a_i l^{2i-1} (l-\frac{2ki}{\nu})+ \sum_{i=1}^N a_i\sum_{j=1}^{i}C_{2i}^{2j}l^{2i-2j}c_j (4kt)^j \\
& & - \frac{2k}{\nu} \sum_{i=2}^N a_i\sum_{j=1}^{i-1}C_{2i}^{2j+1}l^{2i-2j-1}c_{j+1} (4kt)^{j} \\
& = & a_0 + \sum_{i=1}^N a_i l^{2i-1} (l-\frac{2ki}{\nu})  + \sum_{j=1}^N \left[a_j c_j+ \sum_{i=j+1}^{n} a_i  l^{2i-2j-1}\left (lC_{2i}^{2j}c_j -\frac{2k}{\nu}C_{2i}^{2j+1}c_{j+1} \right )\right](4kt)^j
\end{eqnarray*}
We see from above that $g_1(t)$ is a polynomial of $t$ thus if there is a selection of $a_i, i=1,2,\cdots$ such that
$$g_1(t)=T_1(t)=T_0(t)-T_0(0)+\frac{k}{\nu}\mu'(l)+\mu(l)$$
then the boundary conditions in \eqnref{eq:103} are satisfied.
Since $T_0(t)$ is analytic we suppose $T_0(t)= \sum_{i=0}^N d_i t^i$, then
$T_1(t)=\frac{k}{\nu}\mu'(l)+\mu(l)+\sum_{i=1}^N d_i t^i$ and letting $\tilde{d}_0:=\frac{k}{\nu}\mu'(l)+\mu(l)$
the coefficients of the function $\mu(x)$ can be calculated by solving the linear equation \eqnref{matrix}.
We see from (\ref{matrix}) that the matrix is an upper triangle matrix and it is nonsingular since the elements in
the diagonal are nonzero. Furthermore, if $l=\frac{2k}{\nu}$, the matrix is even more simple (the elements in the
secondary diagonal of the matrix are zero). Thus for any analytic function $T_0(t)$ we can construct a function
$\mu(x)$ such that $g_1(t)=T_1(t)$.

Suppose we have solved problem (\ref{eq:103}) and the solution is $u_1(x,t)$,
\begin{eqnarray}
u_1(x,t) & = & \int_{-\infty}^{\infty} \frac{1}{\sqrt{4k\pi t}}e^{-\frac{(x-y)^2}{4kt}} \mu(y) dy  =  \sum_{i=0}^N a_i\frac{1}{\sqrt{4k\pi t}}\int_{-\infty}^{\infty} e^{-\frac{(x-y)^2}{4kt}} y^{2i} dy \nonumber\\
& = & \sum_{i=0}^N a_i\frac{1}{\sqrt{4k\pi t}} \sum_{j=0}^i C_{2i}^{2j} x^{2i-2j} \int_{-\infty}^{\infty} y^{2j}e^{-\frac{y^2}{4kt}}dy = \sum_{i=0}^N a_i\sum_{j=0}^i c_jC_{2i}^{2j} x^{2i-2j}(4kt)^j \label{sol1}
\end{eqnarray}
then the solution to \eqnref{eq:101} is $u(x,t)=u_1(x,t)+u_2(x,t)$, where $u_2(x,t)$ is the solution to
\begin{equation}
\left\{
\begin{array}{ll}
 u_t = ku_{xx} & (x,t) \in \Omega_T, \\
 u(x,0)=\mu_0(x)-\mu(x) & x\in (0,l), \\
 u_x(0,t)=0, \ \ -ku_x(l,t)=\nu (u(l,t)-C_T)  & t \in (0,T],
\end{array}
\right.
\label{eq:106}
\end{equation}
where $C_T=T_0(0)-\frac{k}{\nu}\mu'(l)-\mu(l)$. Denote $\bar{u}(x,t)=u(x,t)-C_T$ then problem (\ref{eq:106}) turns to
\begin{equation}
\left\{
\begin{array}{ll}
 \bar{u}_t = k\bar{u}_{xx} & (x,t) \in \Omega_T, \\
 \bar{u}(x,0)=\mu_0(x)-\mu(x)-C_T & x\in (0,l), \\
 \bar{u}_x(0,t)=0, \ \ -k\bar{u}_x(l,t)=\nu \bar{u}(l,t)  & t \in (0,T].
\end{array}
\right.
\label{eq:a1}
\end{equation}
We can thus use the well-known variational separation method (see, e.g. \cite{Pow:2006}), which is looking for
simple solutions in the form
$$\bar{u}(x,t)=X(x)T(t).$$
By some elementary analysis, we get the solution form to \eqnref{eq:106}
\begin{equation}
\label{eq:108}
u(x,t)= \sum_{n=1}^{\infty} b_n e^{-\sigma_n^2kt}\cos(\sigma_n x) + C_T,
\end{equation}
with $\sigma_n>0, k \sigma_n \tan(\sigma_n l) = \nu$. The choice of $\Gs_n$ is increasing as $n$ increases and
we have $\sigma_n \rightarrow \frac{n\pi}{l}$ as $n\rightarrow \infty$. It is easy to verify that
\eqnref{eq:108} satisfies the boundary conditions. The only problem remaining here is to seek for the constants
$b_n$ such that
\begin{equation}
\label{eq:109}
\sum_{n=1}^{\infty} b_n \cos(\sigma_n x)=\mu_0(x)-\mu(x)-C_T.
\end{equation}
As we know that $\{\cos(\frac{n \pi}{l} x)\}$ is an orthonormal basis in $L^2([0,l])$. We can also prove that
$\{\cos(\sigma_n x)\}$ is orthogonal in $L^2([0,1])$. In fact by using $k \sigma_n \tan(\sigma_n l) = \nu$ we have
\begin{eqnarray*}
\int_0^l \cos(\sigma_n x) \cos(\sigma_m x)dx & = & \frac{1}{\sigma_n} \sin(\sigma_n x)\cos(\sigma_m x)\Big|_0^l
+ \frac{\sigma_m}{\sigma_n} \int_0^l \sin(\sigma_n x) \sin (\sigma_m x) dx \\
& = & \frac{1}{\sigma_n^2}(\sigma_n\sin(\sigma_n l)\cos(\sigma_m l)- \sigma_m \sin (\sigma_m l)\cos(\sigma_n l)) \\
& & + \frac{\sigma_m^2}{\sigma_n^2} \int_0^l \cos(\sigma_n x) \cos (\sigma_m x) dx =
 \frac{\sigma_m^2}{\sigma_n^2} \int_0^l \cos(\sigma_n x) \cos (\sigma_m x) dx.
\end{eqnarray*}
Furthermore, by Robin boundary condition we know that $u(x,0)$
can not be nonzero constant functions. Thus we have
$$b_n= \frac{\int_0^l (\mu_0(x)-\mu(x)-C_T) \cos(\sigma_n x) dx}{\int_0^l \cos^2 (\sigma_n x)dx}
= \frac{2\nu }{\nu l+k \sin^2(\sigma_nl)}\int_0^l (\mu_0(x)-\mu(x)-C_T) \cos(\sigma_n x) dx .$$
We should mention that in order to make the summation in the left hand side of (\ref{eq:109}) converge, some continuous conditions on
$\mu_0(x)$ should be given. However, we do not give
further discussion on this because it can be similarly inherited from the Fourier series theory (see, e.g., \cite{Pow:2006}).

From the analysis above, we conclude that the solution to (\ref{eq:101}) with Neumann-Robin boundary conditions has the
form \eqnref{eq:solnr} if $F(x,t)=0$. The strong solution
to (\ref{eq:101}) with $F(x,t)$ not identically zero can be solved similarly as discussed.
Followed the idea above, we only
need to extend $F(x,t)$ to be an even function in the whole $x$-axis, that is $F(-x,t)=F(x,t)$ then same method can be
used to get the solution of (\ref{eq:101}). In this case, the boundary values are
\begin{eqnarray*}
u(l,t) &  =  & \int_{-\infty}^{\infty} \frac{1}{\sqrt{4k\pi t}}e^{-\frac{(l-y)^2}{4kt}} \mu_0(y) dy
+ \int_{0}^{t}\int_{-\infty}^{\infty} \frac{1}{\sqrt{4k\pi (t-s)}}e^{-\frac{(l-y)^2}{4k(t-s)}} F(y,s) dyds, \\
u_x(l,t) &  =  & -\int_{-\infty}^{\infty} \frac{l-y}{2kt\sqrt{4k\pi t}}e^{-\frac{(l-y)^2}{4kt}} \mu_0(y) dy \\
& & - \int_{0}^{t}\int_{-\infty}^{\infty} \frac{l-y}{2k(t-s)\sqrt{4k\pi (t-s)}}e^{-\frac{(l-y)^2}{4k(t-s)}} F(y,s) dyds.
\end{eqnarray*}
We can actually treat $T_1(t)$ as
\begin{eqnarray}
T_1(t) & = & T_0(t)-\int_{0}^{t}\int_{-\infty}^{\infty} \frac{1}{\sqrt{4k\pi (t-s)}}e^{-\frac{(l-y)^2}{4k(t-s)}} F(y,s) dyds
\nonumber \\
& & - \int_{0}^{t}\int_{-\infty}^{\infty} \frac{l-y}{2\nu (t-s)\sqrt{4k\pi (t-s)}}e^{-\frac{(l-y)^2}{4k(t-s)}} F(y,s) dyds
-d_0+\frac{k}{\nu}\tilde{\mu}'(l)+\tilde{\mu}(l)
\label{eq:111}
\end{eqnarray}
and do some same analysis to get the solution of (\ref{eq:101}).

\subsection{Dirichlet-Robin boundary condition}
To solve the strong solution of (\ref{eq:101}) with Dirichlet-Robin boundary condition, we
shall first consider the solution to
\begin{equation}
\label{eq:201}
\left\{
\begin{array}{ll}
 u_t = ku_{xx} & (x,t) \in \Omega_T, \\
 u(x,0)=\tilde{\mu}(x) & x\in (0,l), \\
 u(0,t)=0, \ \ -k u_x(l,t)=\nu (u(l,t)-T_2(t))& t \in (0,T]. \\
\end{array}
\right.
\end{equation}
 $T_2(t)=T_0(t)-d_0+\frac{k}{\nu}\tilde{\mu}'(l)+\tilde{\mu}(l)$.
In order to satisfy the Dirichlet boundary condition the search of
$\tilde{\mu}(x)$ is now for an odd function of the following form
\begin{equation}
\label{eq:202}
\tilde{\mu}(x)= \sum_{i=0}^N \tilde{a}_i x^{2i+1},
\end{equation}
and by using the Robin boundary we compute $g_2(t)$
\begin{eqnarray*}
g_2(t) & = & \frac{1}{\sqrt{4k \pi t}} \int_{-\infty}^{\infty} \left (1 - \frac{l-y}{2\nu t} \right)
e^{-\frac{(l-y)^2}{4kt}} \sum_{i=0}^N a_i y^{2i+1} dy \\
& = & \frac{1}{\sqrt{\pi}} \sum_{i=0}^N a_i\int_{-\infty}^{\infty} e^{-y^2}  (l-\sqrt{4kt}y)^{2i+1} dy
-  \frac{\sqrt{4kt}}{2\nu t\sqrt{\pi}} \sum_{i=0}^N a_i\int_{-\infty}^{\infty} y e^{-y^2}  (l-\sqrt{4kt}y)^{2i+1} dy \\
& = &  \sum_{j=0}^N \left[\sum_{i=j}^{N} a_i  l^{2i-2j}\left (lC_{2i+1}^{2j}c_j -\frac{2k}{\nu}C_{2i+1}^{2j+1}c_{j+1} \right )\right](4kt)^j.
\end{eqnarray*}
By setting $g_2(t)=T_2(t)$, we thus get the linear equation \eqnref{matrix3}.
The matrix in (\ref{matrix3}) is also an upper triangle matrix with the diagonal elements nonzero if $l\neq \frac{k}{\nu}$.
In the case $l = \frac{k}{\nu}$ we can actually set $\tilde{\mu}(x)=\sum_{i=0}^{n+1}a_ix^{2i+1}$ with $a_0$ arbitrarily given,
since elements of the secondary diagonal of the matrix (\ref{matrix3}) are totally nonzero. In fact, we know if
$u(x,t)$ is the solution of $u_t=u_{xx}$ then $u(x,t)+cx$ with arbitrary constant $c$ is also the solution of
$u_t=u_{xx}$. Furthermore, the boundary conditions are also satisfied.
Then the solution to (\ref{eq:201}) is
\begin{eqnarray}
u_1(x,t) & = & \int_{-\infty}^{\infty} \frac{1}{\sqrt{4k\pi t}}e^{-\frac{(x-y)^2}{4kt}} \mu(y) dy =  \sum_{i=0}^n a_i\frac{1}{\sqrt{4k\pi t}}\int_{-\infty}^{\infty} e^{-\frac{(x-y)^2}{4kt}} y^{2i+1} dy \nonumber\\
& = &\sum_{i=0}^n a_i\sum_{j=0}^i c_jC_{2i+1}^{2j} x^{2i-2j+1}(4kt)^j \label{sol2}
\end{eqnarray}
Next we need further solve
\begin{equation}
\left\{
\begin{array}{ll}
 u_t = ku_{xx} & (x,t) \in \Omega_T, \\
 u(x,0)=\mu_0(x)-\tilde{\mu}(x) & x\in (0,l), \\
 u(0,t)=0, \ \ -ku_x(l,t)=\nu (u(l,t)-\tilde{C}_T)  & t \in (0,T],
\end{array}
\right.
\label{eq:203}
\end{equation}
where $\tilde{C}_T=T_1(0)-\frac{k}{\nu}\tilde{\mu}'(l)-\tilde\mu(l)$. Similarly, we can get the solution to \eqnref{eq:203},
denoted by $u_2(x,t)$
$$u_2(x,t)= \sum_{n=1}^{\infty} \tilde{b}_n e^{-\tilde\sigma_n^2kt}\sin(\tilde\sigma_n x)+\tilde{C}_T,$$
where $\tilde\sigma_n>0$ satisfy $k \tilde\sigma_n  = -\nu\tan(\tilde\sigma_n l)$ and
$$\tilde{b}_n= \frac{\int_0^l (\tilde\mu_0(x)-\tilde\mu(x)-\tilde{C}_T) \sin(\tilde\sigma_n x) dx}{\int_0^l \sin^2 (\tilde\sigma_n x)dx}
= \frac{2\nu }{\nu l+k \cos^2(\tilde\sigma_nl)}\int_0^l (\tilde\mu_0(x)-\tilde\mu(x)-C_T) \sin(\tilde\sigma_n x) dx .$$

\section{Different solution forms}
In this section we study the solution of the following problem
\begin{equation}
\label{eq:501}
\left\{
\begin{array}{ll}
 u_t = ku_{xx} + f(x) , \\
 u(x,0)=\mu_0(x)  \\
 u_x(0,t)=0 , \ \ u_x(1,t)=0.
\end{array}
\right.
\end{equation}
Firstly, by separating variables we have the solution of (\ref{eq:501}) as follows:
\begin{equation}
\label{eq:502}
u(x,t)= \sum_{n=0}^{\infty}\left[e^{-n^2\pi^2kt}a_n \cos(n\pi x) +\frac{1-e^{-n^2\pi^2kt}}{n^2\pi^2k}b_n \cos(n\pi x)\right],
\end{equation}
where for $n=0$, $a_0=\int_0^1 \mu_0(x)dx$, $b_0=\int_0^1 f(x)dx$, and for $n>0$,
$a_n= 2\int_0^1 \mu_0(x) \cos(n\pi x) dx$,  $b_n=2\int_0^1 f(x) \cos(n\pi x) dx$. (\ref{eq:502}) makes no sense
for $n=0$. However we can define
$$\frac{1-e^{-n^2\pi^2kt}}{n^2\pi^2}|_{n=0}:=\lim_{n\rightarrow 0} \frac{1-e^{-n^2\pi^2kt}}{n^2\pi^2k}= t,$$
to make it meaningful at $n=0$. Secondly, by using the extension method in this paper, we can deduce a
different solution form of (\ref{eq:501}) as follows:
\begin{equation}
\label{eq:503}
u(x,t)= \int_{-\infty}^{\infty} \frac{1}{\sqrt{4k\pi t}}e^{-\frac{(x-y)^2}{4kt}} \mu_0(y) dy
+ \int_{0}^{t}\int_{-\infty}^{\infty} \frac{1}{\sqrt{4k\pi s}}e^{-\frac{(x-y)^2}{4ks}} f(y) dyds,
\end{equation}
where $\mu_0(x)$ and $f(x)$ are extended to the following form by Fourier expansion
\begin{eqnarray*}
\mu_0(x) & = &  \int_0^1 \mu_0(y)dy + 2\sum_{n=1}^{\infty} \cos(n\pi x) \int_0^1 \mu_0(y) \cos(n\pi y) dy, \\
f(x) & = & \int_0^1 f(y)dy + 2\sum_{n=1}^{\infty} \cos(n\pi x) \int_0^1 f(y) \cos(n\pi y) dy.
\end{eqnarray*}
The boundary conditions are satisfied automatically by this kind of extension. The solutions (\ref{eq:503}) and
(\ref{eq:502}) have totally different forms while the solution of (\ref{eq:501}) is unique. Thus the two different
solution forms must be the same function. To explain this in detail, we first verify that the values are the same at
$x=0$. By (\ref{eq:502}) we obtain
$$u(0,t)=\sum_{n=0}^{\infty}\left[e^{-n^2\pi^2kt}a_n +
\frac{1-e^{-n^2\pi^2kt}}{n^2\pi^2k}b_n\right].$$
On the other hand by (\ref{eq:503}) we have
\begin{eqnarray*}
u(0,t) & = & \int_{-\infty}^{\infty} \frac{1}{\sqrt{4k\pi t}}e^{-\frac{y^2}{4kt}}
 \sum_{n=0}^{\infty} a_n\cos(n\pi y) dy  + \int_{0}^{t}\int_{-\infty}^{\infty} \frac{1}{\sqrt{4k\pi s}}e^{-\frac{y^2}{4ks}}
\sum_{n=0}^{\infty} b_n\cos(n\pi y) dyds \\
& = & \sum_{n=0}^{\infty} a_n\int_{-\infty}^{\infty} \cos(n\pi y)\frac{1}{\sqrt{4k\pi t}}e^{-\frac{y^2}{4kt}} dy
 +\sum_{n=0}^{\infty} b_n
\int_{0}^{t}\int_{-\infty}^{\infty} \cos(n\pi y)\frac{1}{\sqrt{4k\pi s}}e^{-\frac{y^2}{4ks}}dyds.
\end{eqnarray*}
Comparing the two forms above, it is enough to verify the following
\begin{lemma}
There holds the following
\begin{equation}
\label{eq:504}
\int_{-\infty}^{\infty} cos(n\pi y)\frac{1}{\sqrt{4k\pi t}}e^{-\frac{y^2}{4kt}} dy = e^{-n^2\pi^2kt}.
\end{equation}
\end{lemma}
\pf  By Taylor expansion of $\cos(n\pi y)$ we have
\begin{eqnarray*}
\int_{-\infty}^{\infty} \cos(n\pi y)\frac{1}{\sqrt{4k\pi t}}e^{-\frac{y^2}{4kt}} dy
& = & \int_{-\infty}^{\infty} \sum_{i=0}^{\infty}(-1)^i \frac{(n\pi y)^{2i}}{(2i)!}
\frac{1}{\sqrt{4k\pi t}}e^{-\frac{y^2}{4kt}} dy \\
& = & \sum_{i=0}^{\infty}\frac{(-1)^i (n\pi)^{2i}}{(2i)!} \int_{-\infty}^{\infty} y^{2i}
\frac{1}{\sqrt{4k\pi t}}e^{-\frac{y^2}{4kt}} dy \\
& = & \sum_{i=0}^{\infty}\frac{(-1)^i (n\pi)^{2i}}{(2i)!} (4kt)^i\prod_{j=1}^i \frac{2j-1}{2} \\
& = & \sum_{i=0}^{\infty}\frac{(-1)^i (n\pi)^{2i}}{i!}(kt)^i= e^{-n^2\pi^2kt}.
\end{eqnarray*}
Next, we obtain
\begin{eqnarray*}
& & \int_{-\infty}^{\infty} \cos(n\pi y)\frac{1}{\sqrt{4k\pi t}}e^{-\frac{(x-y)^2}{4kt}} dy \\
 & =& \int_{-\infty}^{\infty} \cos(n\pi (x-y))\frac{1}{\sqrt{4k\pi t}}e^{-\frac{y^2}{4kt}} dy \\
& = &  \int_{-\infty}^{\infty} \left[\cos(n\pi x) \cos(n\pi y)+\sin(n\pi x) \sin(n\pi y)\right]
\frac{1}{\sqrt{4k\pi t}}e^{-\frac{y^2}{4kt}} dy\\
& = & \cos(n\pi x) \int_{-\infty}^{\infty} \cos(n\pi y)\frac{1}{\sqrt{4k\pi t}}e^{-\frac{y^2}{4kt}} dy\\
& = & \cos(n\pi x) e^{-n^2\pi^2kt}.
\end{eqnarray*}
\qed

This Lemma verifies that (\ref{eq:502}) and (\ref{eq:503})
are the same solution to \eqnref{eq:501} but with two different solution forms .

\section{Some Examples}
In this section we shall present some concrete examples to show how
to get the strong solution to \eqnref{eq:101}. We set $l=1$, $k=1/4$ and $\nu=1/2$,
then the problem (\ref{eq:101}) becomes
\begin{equation}
\label{eq:401}
\left\{
\begin{array}{ll}
 u_t = \frac{1}{4}u_{xx} + F(x,t) & (x,t) \in \Omega_T, \\
 u(x,0)=\mu_0(x) & x\in (0,l), \\
 u_x(0,t)=0, \ \ -u_x(1,t)=2(u(1,t)-T_0(t))  & t \in (0,T],
\end{array}
\right.
\end{equation}
and the linear equation in (\ref{matrix}) turns to
\begin{equation}
\label{matrix1}
\begin{bmatrix}
1 & 0 & \cdots & 1-n \\
0 & c_1 & \cdots & C_{2n}^{2}c_1-C_{2n}^3c_2\\
\vdots & \vdots & \ddots & \vdots \\
0 & \cdots & c_{n-1} & C_{2n}^{2n-2}c_{n-1}-C_{2n}^{2n-1}c_n \\
0 & 0 & \cdots & c_n
\end{bmatrix}
\begin{bmatrix}
a_0 \\
a_1 \\
\vdots \\
a_{n-1} \\
a_n \\
\end{bmatrix}
=
\begin{bmatrix}
d \\
d_1 \\
\vdots \\
d_{n-1} \\
d_n
\end{bmatrix}.
\end{equation}
Since $c_n=\frac{2n-1}{2}c_{n-1}$ we get $C_{2n}^{2n-2}c_{n-1}-C_{2n}^{2n-1}c_n=0$.
For $n=4$, the matrix appearing in the linear equation \eqnref{matrix1} reads
$$
\begin{bmatrix}
1 & 0 & -1 & -2 & -3\\
0 & 0.5 & 0 & -7.5 & -28\\
0 & 0 & 0.75 & 0 & -52.5\\
0 & 0 & 0 & 1.875 & 0 \\
0 & 0 & 0 & 0 & 6.5625
\end{bmatrix}
$$
In what follows, we show how to get the strong solution in four different choice of $F(x,t)$, $T_0(t)$ and $\mu_0(t)$.
\begin{enumerate}
\item[\textbf{Ex1:}] $F(x,t)=2t+1$, $T_0(t)=t^2+t$. \\
From (\ref{eq:111}), by unitary extension, we obtain $T_1(t)=0$ and by (\ref{matrix1})
we have $a_0=a_1=a_2=0$. Thus $\mu(x)=0$. In fact, in this case, problem (\ref{eq:101}) is
$$
\left\{
\begin{array}{ll}
 u_t = \frac{1}{4}u_{xx} + 2t+1 & (x,t) \in [0,1]\times[0,T], \\
 u(x,0)=\mu_0(x) & x\in (0,1), \\
 u_x(0,t)=0, \ \ -u_x(1,t)=2(u(1,t)-t^2-t)  & t \in (0,T],
\end{array}
\right. $$
which can be solved by first setting $\tilde{u}(x,t)=u(x,t)-t^2-t$ and then solving the homogeneous problem
$$
\left\{
\begin{array}{ll}
 \tilde{u}_t = \frac{1}{4}\tilde{u}_{xx} & (x,t) \in [0,1]\times[0,T], \\
 \tilde{u}(x,0)=\mu_0(x) & x\in (0,1), \\
 \tilde{u}_x(0,t)=0, \ \ -\tilde{u}_x(l,t)=2\tilde{u}(l,t)  & t \in (0,T].
\end{array}
\right. $$
\item[\textbf{Ex2:}] $F(x,t)=3t^2+2t$, $T_0(t)=t^3+t^2+t+5$, $\mu_0(x)=2x^2+1$. \\
From (\ref{eq:111}), by unitary extension, we have $T_1(t)=t+5$. Thus $a_0=5$, $a_1=2$. $\mu(x)=2x^2+5$.
By (\ref{eq:001}), (\ref{sol1}) and (\ref{eq:111}), the solution
$u_1(x,t)=a_0+ a_1(x^2+\frac{1}{2}t)+t^3+t^2=2x^2+t^3+t^2+t+5$ is the solution of
$$
\left\{
\begin{array}{ll}
 u_t = \frac{1}{4}u_{xx} + 3t^2+2t & (x,t) \in [0,1]\times[0,T], \\
 u(x,0)= 2x^2+5 & x\in (0,1), \\
 u_x(0,t)=0, \ \ -u_x(1,t)=2(u(1,t)-t^3-t^2-t-9)  & t \in (0,T],
\end{array}
\right. $$
And we need further solve
$$
\left\{
\begin{array}{ll}
 u_t = \frac{1}{4}u_{xx} & (x,t) \in [0,1]\times[0,T], \\
 u(x,0)= -4 & x\in (0,1), \\
 u_x(0,t)=0, \ \ -u_x(1,t)=2(u(1,t)+4)  & t \in (0,T],
\end{array}
\right. $$
to get the solution $u_2(x,t)=-4$.
Finally we have $u(x,t)=u_1(x,t)+u_2(x,t)=2x^2+t^3+t^2+t+1$.

On the other hand, if we use the zero extension outside [-1,1] then $T_1(t)=(1-erf(1))(t^3+t^2)+t+5$.
Thus $a_0=5.3775$, $a_1=3.2584$, $a_2=0.20973$, $a_3=0.083893$.
$\mu(x)= 0.083893 x^6+ 0.20973 x^4 + 3.2584 x^2 + 5.3775$ and so $\mu(1)=8.929523$, $\mu'(1)=7.859078$.
By (\ref{sol1}) we know
\begin{eqnarray*}
u_1(x,t) & = & 0.083893 x^6+ 0.20973 x^4 + 3.2584 x^2 + 5.3775  \\
& & +(1.6292+0.6292x^2+0.6292x^4)t + 0.9438x^2t^2 + t^3 +t^2
\end{eqnarray*}
is the solution of
$$
\left\{
\begin{array}{ll}
 u_t = \frac{1}{4}u_{xx} + 3t^2+2t & (x,t) \in [0,1]\times[0,T], \\
 u(x,0)= 0.083893 x^6+ 0.20973 x^4 + 3.2584 x^2 + 5.3775 & x\in (0,1), \\
 u_x(0,t)=0, \ \ -u_x(1,t)=2(u(1,t)-t^3-t^2-t-12.859062)  & t \in (0,T],
\end{array}
\right. $$
and we need further solve
$$
\left\{
\begin{array}{ll}
 u_t = \frac{1}{4}u_{xx} & (x,t) \in [0,1]\times[0,T], \\
 u(x,0)=-0.083893 x^6- 0.20973 x^4 - 1.2584 x^2 - 4.3775 & x\in (0,1), \\
 u_x(0,t)=0, \ \ -u_x(1,t)=2(u(1,t)+ 7.859062) & t \in (0,T].
\end{array}
\right. $$
By elementary calculation we know the solution is approximately
\begin{eqnarray*}
u_2(x,t) & = & -0.083893 x^6- 0.20973 x^4 - 1.2584 x^2 - 4.3775  \\
& & -(0.6292+0.6292x^2+0.6292x^4)t - 0.9438x^2t^2.
\end{eqnarray*}
Finally we also get the solution
$$u(x,t)=u_1(x,t)+u_2(x,t)=2x^2+1+t+t^3+t^2.$$
\item[\textbf{Ex3:}] $F(x,t)=x^2(5t+2)$, $T_0(t)=3t+1$, $\mu_0(x)=x^3+3x^2+1$. \\
If we do not change $F(x,t)$ in the whole $x$-axis, then from (\ref{eq:111}) we obtain $T_1(t)=1+3t-\frac{t^2}{2}-\frac{5}{12}t^3$.
Using (\ref{matrix1}) we obtain $a_0=-\frac{1}{9}$, $a_1=\frac{8}{3}$, $a_2=-\frac{2}{3}$, $a_3=-\frac{2}{9}$. So
$\mu(x)=- \frac{2}{9}x^6 - \frac{2}{3}x^4 + \frac{8}{3}x^2 - \frac{1}{9}$ and $\mu(1)=\frac{5}{3}$, $\mu'(1)=\frac{4}{3}$.
By (\ref{sol1}) we have
\begin{eqnarray*}
u_1(x,t) & = & - \frac{2}{9}x^6- \frac{2}{3}x^4 +\frac{8}{3}x^2 - \frac{1}{9} - t(\frac{5}{3}x^4 - \frac{4}{3})
\end{eqnarray*}
be the solution of
$$
\left\{
\begin{array}{ll}
 u_t = \frac{1}{4}u_{xx} + x^2(5t+2) & (x,t) \in [0,1]\times[0,T], \\
 u(x,0)= - \frac{2}{9}x^6 - \frac{2}{3}x^4 + \frac{8}{3}x^2 - \frac{1}{9} & x\in (0,1), \\
 u_x(0,t)=0, \ \ -u_x(1,t)=2(u(1,t)-3t-\frac{7}{3})  & t \in (0,T],
\end{array}
\right. $$
and we need further solve
$$
\left\{
\begin{array}{ll}
 u_t = \frac{1}{4}u_{xx} & (x,t) \in [0,1]\times[0,T], \\
 u(x,0)=\frac{2}{9}x^6 +  \frac{2}{3}x^4 +3x^3+ \frac{1}{3}x^2 + \frac{10}{9} & x\in (0,1), \\
 u_x(0,t)=0, \ \ -u_x(1,t)=2(u(1,t)+\frac{4}{3}) & t \in (0,T].
\end{array}
\right. $$
\end{enumerate}

\section{Conclusions}
We have shown the existence and uniqueness of the strong solution (or the classical solution) for linear
parabolic equation with mixed boundary conditions. Although we only analyzed the {\em Neumann-Robin} boundary problem
and the {\em Dirichlet-Robin} boudnary problem, the other boundary condition problems can also be solved similarly.
Furthermore, the extension method used in this paper can also be easily implemented as a numerical method for solving
the linear parabolic equation and the efficiency of the method is quite competitive. Even if $T_1(t)$ is not infinitely
differentiable in $[0,l]$ the numerical solution can also be quite accurate if $T_1(t)$ can be approximated well with
certain degrees of polynomials. Future works would focus on the numerical comparison between the method in this paper
and other numerical methods such as finite difference and finite element method. The method can also be extended to two
dimensional or higher dimensional linear parabolic equations in the future.

\end{document}